\documentclass{article}%
\usepackage{amsmath}%
\usepackage{amsfonts}%
\usepackage{amssymb}%
\usepackage{graphicx}

\begin{document}

\title{A Generalization of Fourier Series occurring in Atomic Theory}
\author{Bernard J. Laurenzi\\Department of Chemistry\\The State University of New York at Albany}
\date{June, 29, 2019}
\maketitle

\begin{abstract}
A number of the Fourier series which occur in the theory of the semi-classical
atom due to Englert and Schwinger are generalized and presented.

\end{abstract}

\section{Some Generalized Fourier Series}

In the work shown below we have set the goal of presenting certain Fourier
series in terms of closed form expressions which contain the special functions
commonly found in symbolic computational software such as Maple and
Mathematica. In advance of the derivation of the sums to be obtained below it
will prove useful to introduce the functions $\{z\}$ and $<z>$ i.e.
\[
\{z\}=frac(z)=z-\lfloor z\rfloor\,,\hspace{0.25in}0\leq\{z\}<1,
\]%

\[
<z>\,=z-\lfloor z+1/2\rfloor\,,\hspace{0.25in}-\frac{1}{2}\leq\;<z>\;<\frac{1}%
{2},
\]
where $\{z\}$ is the fractional part of $z$ and $\lfloor z\rfloor$ is the
floor function. We note that%
\[
\{z\}=\,<z-\frac{1}{2}>+\frac{1}{2}.
\]

\subsection{Englert alternating sums}

The alternating sum%
\[
\sum_{k=1}^{\infty}\frac{(-1)^{k}\sin(2\pi kz)}{\pi k}=-<z>,
\]
has been given by Englert, \textit{et al}. \cite{Englert}. \ Using that
expression, the sums given below have been obtained by repeated integration of
that series or obtained from the literature \cite{Ober}, \cite{Hansen} and
represent an extension of Englert's work. \ We have in the first instances
sums which are expressible as polynomials which contain the function $<z>$
i.e.
\begin{align*}
S_{0}(z)  & =\sum_{k=1}^{\infty}\frac{(-1)^{k}\sin(2\pi kz)}{\pi k}=-<z>,\\
C_{1}(z)  & =\sum_{k=1}^{\infty}\frac{(-1)^{k}\cos(2\pi kz)}{(\pi k)^{2}%
}=<z>^{2}-\frac{1}{12},\\
S_{1}(z)  & =\sum_{k=1}^{\infty}\frac{(-1)^{k}\sin(2\pi kz)}{(\pi k)^{3}%
}=\frac{2}{3}<z>[<z>^{2}-\frac{1}{4}],\\
C_{2}(z)  & =\sum_{k=1}^{\infty}\frac{(-1)^{k}\cos(2\pi kz)}{(\pi k)^{4}%
}=-\frac{1}{3}<z>^{2}[<z>^{2}-\frac{1}{2}]-\frac{7}{720},\\
S_{2}(z)  & =\sum_{k=1}^{\infty}\frac{(-1)^{k}\sin(2\pi kz)}{(\pi k)^{5}%
}=-\frac{1}{3}<z>[\frac{2}{5}<z>^{4}-\frac{1}{3}<z>^{2}+\frac{7}{120}],\\
C_{3}(z)  & =\sum_{k=1}^{\infty}\frac{(-1)^{k}\cos(2\pi kz)}{(\pi k)^{6}%
}=\frac{1}{9}<z>^{2}[\frac{2}{5}<z>^{4}-\frac{1}{2}<z>^{2}+\frac{7}%
{40}]-\frac{31}{30240}.
\end{align*}
General expressions for these sum can be obtained and are presented in the
sequel. \ Important in what follows we note the integral relations%
\begin{align}
\int_{0}^{z}  & <\zeta>^{2n+1}d\zeta=\frac{1}{2n+2}<z>^{2n+2},\label{eq1}\\
\int_{0}^{z}  & <\zeta>^{2n}d\zeta=\frac{1}{2n+1}[<z>^{2n+1}+\frac{1}{2^{2n}%
}(z-<z>)].\label{eq2}%
\end{align}
Defining the sums $\mathcal{S}_{n}(z)$ and $\mathcal{C}_{n}(z)$ whose summands
are \textit{even functions} of the index $k\ $by%
\begin{align*}
\mathcal{S}_{n}(z)  & =\sum_{k=1}^{\infty}\frac{(-1)^{k}\sin(2\pi kz)}{(\pi
k)^{2n+1}},\\
\mathcal{C}_{n}(z)  & =\sum_{k=1}^{\infty}\frac{(-1)^{k}\cos(2\pi kz)}{(\pi
k)^{2n}},\hspace{0.25in}n>0,
\end{align*}
we get on integration of $\mathcal{C}_{n}(z)$ and $\mathcal{S}_{n}(z)$
\begin{align}
\mathcal{S}_{n}(z)  & =2\int_{0}^{z}\mathcal{C}_{n}(\zeta)d\zeta
,\label{eqq3}\\
\mathcal{C}_{n+1}(z)  & =-2\int_{0}^{z}\mathcal{S}_{n}(\zeta)d\zeta
-\frac{(2^{2n+1}-1)}{(2n+2)!}|B_{2n+2}|,\label{eq4}%
\end{align}
where the $B_{r}$ are the Bernoulli numbers \cite{Bernoulli}. \ Using
expressions for the initial quantities i.e. \
\begin{align*}
\mathcal{S}_{0}(z)  & =-<z>,\\
C_{1}(z)  & =<z>^{2}-\frac{1}{12},
\end{align*}
all of the higher sums can step-by-step be computed. \ 

It is also possible to obtain general expressions for these sums. \ If \ we
write%
\[
\mathcal{C}_{n}(z)=\sum_{i=0}^{n}c_{i}(n)<z>^{2i},
\]
together with the use of (3) and (4) to obtain the `recursion' relation
\begin{equation}
\mathcal{C}_{n}(z)=-4\int_{0}^{z}\int_{0}^{\zeta}\mathcal{C}_{n-1}%
(\zeta^{\prime})d\zeta^{\prime}d\zeta-\frac{(2^{2n+1}-1)}{(2n+2)!}%
|B_{2n+2}|,\label{eq5}%
\end{equation}
we can show that the coefficients $c_{i}(n)$ contained in the polynomial
expressions for $\mathcal{C}_{n}(z)$ are related to the coefficients
$c_{i}(n-1)$ of the previous polynomial $\mathcal{C}_{n-1}(z)$. \ We get%
\begin{align*}
c_{0}(n)  & =-\frac{(2^{2n+1}-1)}{(2n+2)!}|B_{2n+2}|,\\
c_{1}(n)  & =\sum_{i=0}^{n-1}\frac{c_{i}(n-1)}{2^{2i-1}(2i+1)},\\
c_{i}(n)  & =-\frac{2\,c_{i-1}(n-1)}{i(2i-1)},\hspace{0.25in}i>1,
\end{align*}
and%
\[
\sum_{i=0}^{n-1}\frac{c_{i}(n-1)}{2^{2i}(2i+1)}=0.
\]
Finally we have%
\[
\mathcal{C}_{n}(z)=\left[  \sum_{i=0}^{n-1}\frac{c_{i}(n-1)}{2^{2i-1}%
(2i+1)}\right]  <z>^{2}+2\sum_{i=2}^{n}\frac{\,c_{i-1}(n-1)}{i(2i-1)}%
<z>^{2\,i}-\frac{(2^{2n+1}-1)}{(2n+2)!}|B_{2n+2}|.
\]
In a similar way we also have%
\[
\mathcal{S}_{n}(z)=2\sum_{i=1}^{n-1}\frac{c_{i}(n)}{(2i+1)}[<z>^{2i}%
-\frac{1}{4^{i}}].
\]

In the second instance, the remaining alternating sums \ which are \textit{odd
functions} of the index $k$ will be seen to contain expressions which are more
complicated . \ Writing the sums as \
\begin{align*}
\mathcal{S}_{n}^{\prime}(z)  & =\sum_{k=1}^{\infty}\frac{(-1)^{k}\sin(2\pi
kz)}{(\pi k)^{2n}},\\
\mathcal{C}_{n}^{\prime}(z)  & =\sum_{k=1}^{\infty}\frac{(-1)^{k}\cos(2\pi
kz)}{(\pi k)^{2n+1}},
\end{align*}
we have for $\mathcal{S}_{0}^{\prime}(z)$ \cite{Hardy}%
\[
\mathcal{S}_{0}^{\prime}(z)=-\frac{1}{2}\tan(\pi z),
\]
and for $\mathcal{C}_{0}^{\prime}(z)$ \cite{Hansen1}%
\[
\mathcal{C}_{0}^{\prime}(z)=-\frac{1}{\pi}\ln|2\cos(\pi z)|,
\]
then (see the Appendix)
\begin{align*}
\mathcal{S}_{n}^{\prime}(z)  & =\frac{1}{\pi^{2n}}\operatorname{Im}%
[Li_{2n}(-\exp(-2\pi iz))],\text{ \ \ }n>0,\\
\mathcal{C}_{n}^{\prime}(z)  & =\frac{1}{\pi^{2n+1}}\operatorname{Re}%
[Li_{2n+1}(-\exp(-2\pi iz))],
\end{align*}
where $Li_{a}(z)$ is the polylogarithm function of order $a$ \cite{Lewin}%
\[
Li_{a}(z)=\sum_{k=1}^{\infty}\frac{z^{k}}{k^{a}},
\]%
\[
Li_{1}(z)=-\ln(1-z).
\]
The latter quantities being one of the standard functions contained in Maple
and Mathematica.

\subsection{Non-alternating Englert sums}

The non-alternating sum%
\[
\widetilde{\mathcal{S}}_{0}(z)=\sum_{k=1}^{\infty}\frac{\sin(2\pi kz)}{\pi
k}=\frac{1}{2}-\{z\}=-<z-1/2>,
\]
is due to Titchmarch \cite{Titchmarch}. \ Then%
\[
\widetilde{\mathcal{C}}_{1}(z)=\sum_{k=1}^{\infty}\frac{\cos(2\pi kz)}{(\pi
k)^{2}}=<z-1/2>^{2}-\frac{1}{12},
\]
More generally we define the families of non-alternating \textit{even} and
\textit{odd} sums by \
\begin{align*}
\widetilde{\mathcal{S}}_{n}(z)  & =\sum_{k=1}^{\infty}\frac{\sin(2\pi
kz)}{(\pi k)^{2n+1}},\\
\widetilde{\mathcal{C}}_{n}(z)  & =\sum_{k=1}^{\infty}\frac{\cos(2\pi
kz)}{(\pi k)^{2n}},
\end{align*}
with \ \cite{Hardy2}
\[
\widetilde{\mathcal{C}}_{0}(z)=-\frac{1}{2},
\]
and%
\begin{align*}
\widetilde{\mathcal{S}}_{n}^{\prime}(z)  & =\sum_{k=1}^{\infty}\frac{\sin(2\pi
kz)}{(\pi k)^{2n}},\\
\widetilde{\mathcal{C}}_{n}^{\prime}(z)  & =\sum_{k=1}^{\infty}\frac{\cos(2\pi
kz)}{(\pi k)^{2n+1}},
\end{align*}
respectively. \ We also note the general relations of these sums to the
corresponding alternating \textit{even} and \textit{odd }sums given above
i.e.
\begin{align*}
\widetilde{\mathcal{S}}_{n}(z)  & =\mathcal{S}_{n}(z-1/2),\\
\widetilde{\mathcal{C}}_{n}(z)  & =\mathcal{C}_{n}(z-1/2),\\
\widetilde{\mathcal{S}}_{n}^{\prime}(z)  & =\mathcal{S}_{n}^{\prime}(z-1/2),\\
\widetilde{\mathcal{C}}_{n}^{\prime}(z)  & =\mathcal{C}_{n}^{\prime}(z-1/2).
\end{align*}
Those relations follow directly from the definitions of the `hatted' functions
or the integrated expression of $\widetilde{\mathcal{S}}_{n}(z)$ and
$\widetilde{\mathcal{C}}_{n}(z).$ Here we get%

\[
\widetilde{\mathcal{C}}_{n}(z)=-2\int_{0}^{z}\widetilde{\mathcal{S}}_{n}%
(\zeta)d\zeta+\frac{2^{2n+1}}{(2n+2)!}|B_{2n+2}|,
\]%
\[
\widetilde{\mathcal{S}}_{n}(z)=2\int_{0}^{z}\widetilde{\mathcal{C}}%
_{n}(z)d\zeta.
\]
In this instance we note the integral relations%
\begin{align}
\int_{0}^{z}  & <\zeta-1/2>^{2n+1}d\zeta=\frac{1}{2n+2}\left[  <z-1/2>^{2n+2}%
-\frac{1}{2^{2n+2}}\right]  ,\label{eq6}\\
\int_{0}^{z}  & <\zeta-1/2>^{2n}d\zeta=\frac{1}{2n+1}[<z-1/2>^{2n+1}%
+\frac{1}{2^{2n}}(z-<z>)],\hspace{0.25in}n>0.\label{eq7}%
\end{align}
Using (6), (7) the sums listed below are seen to be polynomials in the
functions $<z-1/2>$. \ That is to say, the non-alternating sums are just the
sums $\mathcal{S}_{n}(z),$ $\mathcal{C}_{n}(z)$ and $\mathcal{S}_{n}^{\prime
}(z),$ $\mathcal{C}_{n}^{\prime}(z)$ which contain the shifted variable i.e.
$z\rightarrow z-1/2.$ \ We have
\begin{align*}
\widetilde{\mathcal{S}}_{0}(z)  & =\sum_{k=1}^{\infty}\frac{\sin(2\pi kz)}{\pi
k}=\frac{1}{2}-\{z\}=-<z-1/2>,\\
\widetilde{\mathcal{C}}_{1}(z)  & =\sum_{k=1}^{\infty}\frac{\cos(2\pi
kz)}{(\pi k)^{2}}=<z-\tfrac{1}{2}>^{2}-\frac{1}{12},\\
\widetilde{\mathcal{S}}_{1}(z)  & =\sum_{k=1}^{\infty}\frac{\sin(2\pi
kz)}{(\pi k)^{3}}=\frac{2}{3}<z-\tfrac{1}{2}>[<z-\tfrac{1}{2}>^{2}-\frac{1}%
{4}],\\
\widetilde{\mathcal{C}}_{2}(z)  & =\sum_{k=1}^{\infty}\frac{\cos(2\pi
kz)}{(\pi k)^{4}}=-\frac{1}{3}<z-1/2>^{2}[<z-1/2>^{2}-\frac{1}{2}%
]-\frac{7}{720},
\end{align*}
and%

\begin{align*}
\widetilde{\mathcal{S}}_{0}^{\prime}(z)  & =\frac{1}{2}\cot(\pi z),\\
\widetilde{\mathcal{C}}_{0}^{\prime}(z)  & =-\frac{1}{\pi}\ln|2\sin(\pi z)|,\\
\widetilde{\mathcal{S}}_{n}^{\prime}(z)  & =\frac{1}{\pi^{2n}}%
\operatorname{Im}[Li_{2n}(\exp(2\pi iz))]\hspace{0.25in}n>0,\\
\widetilde{\mathcal{C}}_{n}^{\prime}(z)  & =\frac{1}{\pi^{2n+1}}%
\operatorname{Re}[Li_{2n+1}(\exp(2\pi iz))]\hspace{0.25in}n>0.
\end{align*}

\subsection{Series with summands containing arguments $2k+1$}

We write the alternating and non-alternating sums which are \textit{even} and
\textit{odd} functions in the argument $(2k+1)$ as%
\begin{align*}
\mathbf{S}_{n}(z)  & =\sum_{k=0}^{\infty}\frac{(-1)^{k}\sin(2\pi
\lbrack2k+1]z)}{[\pi(2k+1)]^{2n+1}},\\
\mathbf{C}_{n}(z)  & =\sum_{k=0}^{\infty}\frac{(-1)^{k}\cos(2\pi
\lbrack2k+1]z)}{[\pi(2k+1)]^{2n}},\\
& \\
\mathbf{S}_{n}^{\prime}(z)  & =\sum_{k=0}^{\infty}\frac{(-1)^{k}\sin
(2\pi\lbrack2k+1]z)}{[\pi(2k+1)]^{2n}},\\
\mathbf{C}_{n}^{\prime}(z)  & =\sum_{k=0}^{\infty}\frac{(-1)^{k}\cos
(2\pi\lbrack2k+1]z)}{[\pi(2k+1)]^{2n+1}},
\end{align*}
and%
\begin{align*}
\widetilde{\mathbf{S}}_{n}(z)  & =\sum_{k=0}^{\infty}\frac{\sin(2\pi
\lbrack2k+1]z)}{[\pi(2k+1)]^{2n+1}},\\
\widetilde{\mathbf{C}}_{n}(z)  & =\sum_{k=0}^{\infty}\frac{\cos(2\pi
\lbrack2k+1]z)}{[\pi(2k+1)]^{2n}},\\
& \\
\widetilde{\mathbf{S}}_{n}^{\prime}(z)  & =\sum_{k=0}^{\infty}\frac{\sin
(2\pi\lbrack2k+1]z)}{[\pi(2k+1)]^{2n}},\\
\widetilde{\mathbf{C}}_{n}^{\prime}(z)  & =\sum_{k=0}^{\infty}\frac{\cos
(2\pi\lbrack2k+1]z)}{[\pi(2k+1)]^{2n+1}},
\end{align*}
respectively. \ 

We find that the alternating and non-alternating sums are interrelated i.e.%
\begin{align}
\widetilde{\mathbf{S}}_{n}(z)  & =\mathbf{C}_{n}^{\prime}(z-1/4)\nonumber\\
\widetilde{\mathbf{C}}_{n}(z)  & =\mathbf{S}_{n}^{\prime}(z+1/4),\nonumber\\
\widetilde{\mathbf{S}}_{n}^{\prime}(z)  & =\mathbf{C}_{n}(z-1/4),\nonumber\\
\widetilde{\mathbf{C}}_{n}^{\prime}(z)  & =\mathbf{S}_{n}(z+1/4).\label{eq88}%
\end{align}
As seen above repeated integrations produce the higher sums. The first few of
the \textit{even} alternating sums being (cf. Appendix)%
\begin{align*}
\mathbf{S}_{0}(z)  & =\frac{1}{2\pi}\ln|\frac{\cos(\pi z)+\sin(\pi z)}%
{\cos(\pi z)-\sin(\pi z)}|=\frac{1}{2\pi}\ln|\tan(\pi z+\pi/4)|,\\
\mathbf{S}_{n}(z)  & =-\frac{1}{2\pi^{2n+1}}\operatorname{Re}[Li_{2n+1}%
(\exp(2\pi i\{z+\tfrac{1}{4}\})-Li_{2n+1}(\exp(2\pi i\{z-\tfrac{1}%
{4}\})],\hspace{0.25in}n>0,\\
\mathbf{C}_{n}(z)  & =\frac{1}{2\pi^{2n}}\operatorname{Im}[Li_{2n}(\exp(2\pi
i\{z+\tfrac{1}{4}\})-Li_{2n}(\exp(2\pi i\{z-\tfrac{1}{4}\})],\hspace
{0.25in}n>0,
\end{align*}
and for the \textit{odd} alternating sums
\begin{align*}
\mathbf{C}_{0}^{\prime}(z)  & =\frac{1}{4}(-1)^{\lfloor2z+1/2\rfloor},\\
\mathbf{S}_{1}^{\prime}(z)  & =\frac{1}{2}[<z+\frac{1}{4}>^{2}-\,<z-\frac{1}%
{4}>^{2}],\\
\mathbf{C}_{1}^{\prime}(z)  & =\frac{1}{3}[<z+\frac{1}{4}>\{\frac{1}%
{4}-\,<z+\frac{1}{4}>^{2}\}-<z-\frac{1}{4}>\{\frac{1}{4}-\,<z-\frac{1}{4}%
>^{2}\}],\\
\mathbf{S}_{2}^{\prime}(z)  & =\frac{1}{6}[<z+\frac{1}{4}>^{2}\{\frac{1}%
{2}-\,<z+\frac{1}{4}>^{2}\}-<z-\frac{1}{4}>^{2}\{\frac{1}{2}-\,<z-\frac{1}%
{4}>^{2}\}].
\end{align*}
The first few of the non-alternating sums can then\ be obtained from (8).

\subsection{Modified Englert sums}

Here we consider sums in which the arguments of the trigonometric functions in
the summands are $2\pi kz$ whereas the arguments in the denominator are
$(2k+1)$ i.e.
\begin{align}
P_{n}(z)  & =\sum_{k=0}^{\infty}\frac{(-1)^{k}\sin(2\pi kz)}{[\pi
(2k+1)]^{2n+1}},\nonumber\\
Q_{n}(z)  & =\sum_{k=0}^{\infty}\frac{(-1)^{k}\cos(2\pi kz)}{[\pi(2k+1)]^{2n}%
},\nonumber\\
& \nonumber\\
P_{n}^{\prime}(z)  & =\sum_{k=0}^{\infty}\frac{(-1)^{k}\sin(2\pi kz)}%
{[\pi(2k+1)]^{2n}},\nonumber\\
Q_{n}^{\prime}(z)  & =\sum_{k=0}^{\infty}\frac{(-1)^{k}\cos(2\pi kz)}%
{[\pi(2k+1)]^{2n+1}},\label{eqqq9}%
\end{align}
and%
\begin{align}
\widetilde{P}_{n}(z)  & =\sum_{k=0}^{\infty}\frac{\sin(2\pi kz)}%
{[\pi(2k+1)]^{2n+1}},\nonumber\\
\widetilde{Q}_{n}(z)  & =\sum_{k=0}^{\infty}\frac{\cos(2\pi kz)}%
{[\pi(2k+1)]^{2n}},\nonumber\\
& \nonumber\\
\widetilde{P}_{n}^{\prime}(z)  & =\sum_{k=0}^{\infty}\frac{\sin(2\pi kz)}%
{[\pi(2k+1)]^{2n}},\nonumber\\
\widetilde{Q}_{n}^{\prime}(z)  & =\sum_{k=0}^{\infty}\frac{\cos(2\pi kz)}%
{[\pi(2k+1)]^{2n+1}}.\label{eq10}%
\end{align}
From (9), and (10) above we get
\begin{align*}
P_{n}(z)  & =\cos(\pi z)\mathbf{S}_{n}(z/2)-\sin(\pi z)\mathbf{C}_{n}^{\prime
}(z/2),\\
Q_{n}(z)  & =\sin(\pi z)\mathbf{S}_{n}^{\prime}(z/2)+\cos(\pi z)\mathbf{C}%
_{n}(z/2),\hspace{0.25in}n>0,\\
P_{n}^{\prime}(z)  & =\cos(\pi z)\mathbf{S}_{n}^{\prime}(z/2)-\sin(\pi
z)\mathbf{C}_{n}(z/2),\hspace{0.25in}n>0,\\
Q_{n}^{\prime}(z)  & =\sin(\pi z)\mathbf{S}_{n}(z/2)+\cos(\pi z)\mathbf{C}%
_{n}^{\prime}(z/2),\\
& \\
\widetilde{P}_{n}(z)  & =\cos(\pi z)\widetilde{\mathbf{S}}_{n}(z/2)-\sin(\pi
z)\widetilde{\mathbf{C}}_{n}^{\prime}(z/2)\\
\widetilde{Q}_{n}(z)  & =\sin(\pi z)\widetilde{\mathbf{S}}_{n}^{\prime
}(z/2)+\cos(\pi z)\widetilde{\mathbf{C}}_{n}(z/2),\hspace{0.25in}n>0,\\
\widetilde{P}_{n}^{\prime}(z)  & =\cos(\pi z)\widetilde{\mathbf{S}}%
_{n}^{\prime}(z/2)-\sin(\pi z)\widetilde{\mathbf{C}}_{n}(z/2),\hspace
{0.25in}n>0,\\
\widetilde{Q}_{n}^{\prime}(z)  & =\sin(\pi z)\widetilde{\mathbf{S}}%
_{n}(z/2)+\cos(\pi z)\widetilde{\mathbf{C}}_{n}^{\prime}(z/2).
\end{align*}
We see that the sums above are not new in that they are related to the sums
$\mathbf{S}_{n}(z),$ $\mathbf{C}_{n}(z),$ $\mathbf{S}_{n}^{\prime}(z),$and
$\mathbf{C}_{n}^{\prime}(z)$ together with $\ \widetilde{\mathbf{S}%
}(z),\widetilde{\mathbf{C}}_{n}(z),\widetilde{\mathbf{S}}_{n}^{\prime}(z),$
and $\widetilde{\mathbf{C}}_{n}^{\prime}(z)$\textbf{. \ }The first few of
these sums are
\begin{align*}
P_{0}(z)  & =-\frac{1}{\pi}\operatorname{Im}\{\exp(\pi iz)\arctan(\exp(-\pi
iz))\},\\
Q_{1}(z)  & =\frac{1}{2\pi^{2}}[\cos(\pi z)\operatorname{Im}\{Li_{2}%
(i\exp(-\pi iz))+Li_{2}(i\exp(\pi iz))\}\\
& +\sin(\pi z\operatorname{Re}\{Li_{2}(i\exp(-\pi iz))-Li_{2}(i\exp(\pi
iz))\}],\\
P_{1}^{\prime}(z)  & =-\frac{1}{2\pi^{2}}[\cos(\pi z)\operatorname{Re}%
\{Li_{2}(i\exp(\pi iz))-Li_{2}(-i\exp(\pi iz))\}\\
& +\sin(\pi z)\operatorname{Im}\{Li_{2}(i\exp(\pi iz))-Li_{2}(-i\exp(\pi
iz))\}],\\
Q_{0}^{\prime}(z)  & =\frac{1}{4}(-1)^{\lfloor z+1/2\rfloor}\cos(\pi
z)-\frac{1}{2\pi}\sin(\pi z)\ln\left|  \frac{\cos(\pi z)}{1-\sin(\pi
z)}\right|  .\\
.  &
\end{align*}
and

\begin{center}%
\begin{align*}
\widetilde{P}_{0}(z)  & =\frac{(-1)^{\lfloor z\rfloor}}{4}\cos(\pi
z)-\frac{1}{2\pi}\sin(\pi z)\ln\left|  \cot(\pi z/2)\right|  ,\\
\widetilde{Q}_{1}(z)  & =\frac{1}{2\pi^{2}}[\cos(\pi z)\operatorname{Re}%
\{Li_{2}(\exp(\pi iz))-Li_{2}(-\exp(\pi iz))\}\\
& +\sin(\pi z)\operatorname{Im}\{Li_{2}(\exp(\pi iz))-Li_{2}(-\exp(\pi
iz))\}],\\
\widetilde{Q}_{0}^{\prime}(z)  & =\frac{1}{4}|\sin(\pi z)|+\frac{1}{2\pi}%
\cos(\pi z)\ln|\cot(\pi z/2)|,\\
\widetilde{P}_{1}^{\prime}(z)  & =\frac{1}{2\pi^{2}}[\cos(\pi
z)\operatorname{Im}\{Li_{2}(\exp(\pi iz))-Li_{2}(-\exp(\pi iz))\}\\
& -\sin(\pi z)\operatorname{Re}\{Li_{2}(\exp(\pi iz))-Li_{2}(-\exp(\pi
iz))\}].\\
&
\end{align*}
\appendix  Appendix
\end{center}

In this appendix we show that in special cases the polylogarith functions
found in the relations above can be written as polynomials in the variables
$<f(z)>$. \ We have%
\[
Li_{a}(\exp(i\theta))=\sum_{k=1}^{\infty}\frac{\cos(k\theta)}{k^{a}}%
+i\sum_{k=1}^{\infty}\frac{\sin(k\theta)}{k^{a}},
\]
from which it follows that%
\begin{align*}
\frac{\operatorname{Re}[Li_{2n}(-\exp(\pm2\pi iz))}{\pi^{2n}}  &
=C_{n}(z)=f_{1,n}(<z>),\\
\frac{\operatorname{Im}[Li_{2n}(-\exp(\pm2\pi iz))}{\pi^{2n}}  & =\pm
S_{n}^{\prime}(z),\\
\frac{\operatorname{Re}[Li_{2n+1}(-\exp(\pm2\pi iz))}{\pi^{2n+1}}  &
=C_{n}^{\prime}(z),\\
\frac{\operatorname{Im}[Li_{2n+1}(-\exp(\pm2\pi iz))}{\pi^{2n+1}}  & =\pm
S_{n}(z)=f_{2,n}(<z>).
\end{align*}
In a similar way we have%
\begin{align*}
\frac{\operatorname{Re}[Li_{2n}(\exp(\pm2\pi iz))}{\pi^{2n}}  & =\widetilde
{C}_{n}(z)=f_{3,n}(<z-1/2>),\\
\frac{\operatorname{Im}[Li_{2n}(\exp(\pm2\pi iz))}{\pi^{2n}}  & =\pm
\widetilde{S}_{n}^{\prime}(z),\\
\frac{\operatorname{Re}[Li_{2n+1}(\exp(\pm2\pi iz))}{\pi^{2n+1}}  &
=\widetilde{C}_{n}^{\prime}(z),\\
\frac{\operatorname{Im}[Li_{2n+1}(\exp(\pm2\pi iz))}{\pi^{2n+1}}  &
=\pm\widetilde{S}_{n}(z)=f_{4,n}(<z-1/2>).
\end{align*}
In the cases where $z$ has been replaced by $z\pm1/4$ the relations are more
complicated i.e.%
\begin{align*}
\frac{\operatorname{Re}[Li_{2n}(\exp(2\pi i[z\pm1/4]))}{\pi^{2n}}  &
=\mp\mathbf{S}_{n}^{\prime}(z)+\frac{1}{2^{2n}}C_{n}(2z)=f_{5,n}(<z>),\\
\frac{\operatorname{Im}[Li_{2n}(\exp(2\pi i[z\pm1/4]))}{\pi^{2n}}  &
=\frac{1}{2^{2n}}S_{n}^{\prime}(2z)\pm\mathbf{C}_{n}(z),\\
\frac{\operatorname{Re}[Li_{2n+1}(\exp(2\pi i[z\pm1/4]))}{\pi^{2n+1}}  &
=\mp\mathbf{S}_{n}(z)+\frac{1}{2^{2n+1}}C_{n}^{\prime}(2z),\\
\frac{\operatorname{Im}[Li_{2n+1}(\exp(2\pi i[z\pm1/4]))}{\pi^{2n+1}}  &
=\frac{1}{2^{2n+1}}S_{n}(2z)\pm\mathbf{C}_{n}^{\prime}(z)=f_{6,n}(<z>).
\end{align*}
The functions $f_{i,n}$ are used\ here to indicate that these Real and
Imaginary parts of the $Li_{a}$ functions can be expressed in terms of
polynomials in $<F(z)>$ \ where $F(z)$ is some linear function of $z$. \ In
the remaining equations above, this does not appear to be possible even in the
cases of infinite series in $<F(z)>$.

\end{document}